\documentclass{amsart}

\usepackage[dvips]{graphics}
\usepackage{amsthm,amssymb,latexsym, indentfirst}
\usepackage{epsfig}
\usepackage{color}
\usepackage{eufrak}
\usepackage[mathscr]{eucal}

\newtheorem*{theorem*}{Theorem}
\newtheorem{proposition}{Proposition}

\newtheorem{lemma}{Lemma}
\newtheorem{fact}{Fact}

\theoremstyle{definition}
\newtheorem{definition}{Definition}

\theoremstyle{remark}
\newtheorem{remark}{Remark}

\newcommand{\PPP}{\mathbb{P}}
\newcommand{\EEE}{\mathbb{E}}

\newcommand{\nid}{\noindent}

\renewcommand{\qedsymbol}{\ensuremath{\blacksquare}}

\newcommand{\B}[1]{\boldsymbol}
\newcommand{\C}[1]{\mathcal{#1}}
\newcommand{\D}[1]{\mathbb{#1}}

\newcommand{\ga}{\alpha}
\newcommand{\gb}{\beta}

\newcommand{\gd}{\delta}
\newcommand{\gep}{\varepsilon}
\newcommand{\gz}{\zeta}

\newcommand{\gl}{\lambda}
\newcommand{\gs}{\sigma}

\newcommand{\gG}{\Gamma}

\newcommand{\ol}[1]{\overline{#1}}
\newcommand{\ul}[1]{\underline{#1}}

\begin{document}

\title[Diffusion in an asymmetric random environment]
{One dimensional diffusion in an asymmetric random environment}
\author{Dimitrios Cheliotis}

\address{
Department of Mathematics \\ Bahen Center for Information
Technology \\ 40 St. George St., 6th floor \\ Toronto, ON M5S 3G3
\\ Canada  }

 \email{dimitris@math.toronto.edu}
 \urladdr{http://www.math.toronto.edu/dimitris/}
\thanks{Research partially supported by an anonymous Stanford Graduate Fellowship and by NSF grant
DMS-0072331  } 
\date{May 27, 2005}

\begin{abstract}
According to a theorem of S. Schumacher, for a diffusion $X$ in an
environment determined by a stable process that belongs to an
appropriate class and has index $a$, it holds that $X_t/(\log
t)^a$ converges in distribution as $t\to\infty$ to a random
variable having an explicit description in terms of the
environment. We compute the density of this random variable in the
case the stable process is spectrally one-sided. This computation
extends a result of H. Kesten and quantifies the bias that the
asymmetry of the environment causes to the behavior of the
diffusion.
\end{abstract}

\maketitle

\section{Introduction} \label{intro}

On the space $\C{W}:=\{f\in \D{R}^{\D{R}}:\text{ $f$ is right
continuous with left limits}\}$ consider the Skorohod topology,
the $\gs$-field of the Borel sets, and $\PPP$ a measure on
$\C{W}$.

Also let $\Omega:=C([0,+\infty))$, and equip it with the
$\gs$-field of Borel sets derived from the topology of uniform
convergence on compact sets. For $w\in\C{W}$, we denote by
$\textbf{P}_w$ the probability measure on $\Omega$ such that
$\{\omega(t):t\ge0\}$ is a diffusion with $\omega(0)=0$ and
generator
\begin{equation}\label{generator}
\frac{1}{2}e^{w(x)}\frac{d}{dx}\left(e^{-w(x)}\frac{d}{dx}\right).
\end{equation}
The construction of such a diffusion is done with a scale and time
transformation from a one-dimensional Brownian motion (see e.g.
\cite{BR}, \cite{SH1}). The diffusion does not explode in finite
time if and only if $k_w(+\infty)=k_w(-\infty)=+\infty$, where
$k_w$ is defined for all $x\in\D{R}$ by $k_w(x):=\int_0^x \int_0^y
e^{w(y)-w(z)}dydz.$ The last statement is Theorem 3 of \cite{SEI}.
We will only consider measures $\PPP$ on $\C{W}$ with the property
\begin{equation} \label{noexplosion}
\PPP\big(k_w(+\infty)=k_w(-\infty)=+\infty\big)=1.
\end{equation}
 For $\PPP$-almost all
$w\in\C{W}$, $\omega$ satisfies the formal SDE
\begin{equation}\label{FSDE}
\begin{array}{rl}
d\omega(t)=&d\beta(t)-\frac{1}{2}w'(\omega(t))\,dt, \\
\omega(0)=&0,
\end{array}
\end{equation}
where $\beta$ is a one-dimensional standard Brownian motion.

Then consider the space $\C{W}\times\Omega$, equip it with the
product $\gs$-field, and take the probability measure defined by
\[
d\C{P}(w,\omega)=d\textbf{P}_w(\omega)\,d\PPP(w).
\]
The marginal of $\C{P}$ in $\Omega$ gives a process which is known
as diffusion in random environment; the environment being the
function $w$.

The following result, concerning a class of environments
$(w_t)_{t\in\D{R}}$, was proved by S. Schumacher (\cite{SC1,
SC2}).
\begin{fact} \label{SchumProp}
Assume \eqref{noexplosion} holds and that there is an $a>0$ such
that the net of processes $\{(w(st)/t^{1/a})_{s\in\D{R}}:t>0\}$
converges to $(U_s)_{s\in\D{R}}$ in the Skorohod topology in
$\C{W}$ as $t\to+\infty$, where $U$ satisfies
\begin{enumerate}
\item [(i)] $U$ is non-degenerate. \item [(ii)] If $[x,y]$ is a
finite interval, and if $\xi:=\inf\{s\ge x:U_s=\inf_{r\in[x,
y]}U_r\}$, then $U$ is continuous at $\xi$. \item[(iii)] If
$[x,y]$ is a finite interval, $U$ attains each of the values
$\inf_{r\in[x, y]}U_r$, $\sup_{r\in[x, y]}U_r$ only once in
$[x,y]$. \item[(iv)] $\PPP(\exists s>0 \text{ such that }
U_s>0)=\PPP(\exists s<0 \text{ such that } U_s>0)=1$.
\end{enumerate}
Then there is a process $b:[0,\infty)\times\C{W}\to\D{R}$ such
that for the formal solution $\omega$ of \eqref{FSDE} it holds
\begin{equation}\label{Schum}
\frac{\omega_t}{(\log t)^a}-b_1(w^{(\log t)})\to 0 \text{ in
$\C{P}$ as $t\to+\infty$,}
\end{equation}
where for $r>0$ we let $w^{(r)}(s)=r^{-1}w(sr^a) \text{ for all
$s\in\D{R}$} .$
\end{fact}

This result shows the dominant effect of the environment, through
the process $b$, on the asymptotic behavior of the diffusion.

In the case where $w$ is a two sided stable process with index
$a\in(1,2]$ and having no positive jumps, we have
$w^{(r)}\overset{\text{law}}{=}w$. Assuming that the assumptions
of the Fact 1 above are satisfied (we will prove this later), we
get $\omega_t/(\log t)^a\Rightarrow b_1(w)$ as $t\to+\infty$. The
main result of this paper is the computation of the density of the
random variable $b_1(w)$.

Let $E_a(z):=\sum_{n=0}^\infty z^n/\Gamma(1+an)$. In the next
section we will introduce two non-negative random variables with
distribution functions $F_u$, $F_d$ respectively, having Laplace
transforms
\begin{align*}
\int_0^\infty e^{-\gl
x}dF_u(x)&=\Big(E_a^{\prime}(\gl)+\frac{a}{a-1}\gl E_a^{\prime\prime}(\gl)\Big)^{-1}\left(\Gamma(a+1)\right)^{-1},\\
\int_0^\infty e^{-\gl
x}dF_d(x)&=\Gamma(a+1)\Big(E_a^{\prime}(\gl)+\frac{a}{a-1}\gl
E_a^{\prime\prime}(\gl)-\frac{a}{a-1}\gl
\frac{\big(E_a^{\prime}(\gl)\big)^2}{E_a(\gl)}\Big),
\end{align*}
for all $\gl>0$. And our main result is the following.

\begin{theorem*} Let $w$ be a L\'evy process with $\EEE(e^{\gl
w_t})=e^{t\gl^a}$ for all $\gl, t\ge0$, where $a$ belongs to
$(1,2]$. The density of $b_1(w)$ is
\begin{equation*}
f_{b_1}(x)=
 \begin{cases}
(a-1)\Gamma(a)\big(1-F_u(-x)\big) &\mbox{ for } x\le0, \\
(a-1)\Gamma(a)\big(1-F_d(x)\big) &\mbox{ for } x>0.
 \end{cases}
\end{equation*}
\end{theorem*}

\begin{remark}\label{biasCh1}
When $a\ne2$, the asymmetry of the environment has a visible
effect on the path of the diffusion. We will show in Section
\ref{hittingtimescomputations} that $\PPP(b_1<0)=\gamma(a)$ for a
strictly decreasing function $\gamma$ of $a$ in $[1,2]$ having
$\gamma(1)=1$ and $\gamma(2)=1/2$. Since, by the theorem of
Schumacher (Fact \ref{SchumProp}),
$\lim_{t\to+\infty}\C{P}(\omega_t<0)=\PPP(b_1<0)$, it follows that
the diffusion is biased towards the left, and the bias increases
as $a$ decreases to 1. Of course, the diffusion is recurrent.
\end{remark}

\begin{remark}
In the case $a=2$, the process $w/\sqrt{2}$ is a standard two
sided Brownian motion and $E_2(z)=\cosh (\sqrt{z})$. The
distribution functions $F_u$, $F_d$ coincide, and their Laplace
transform, appearing above, equals $1/\cosh (\sqrt{\gl})$. Using
Laplace inversion, we recover by the above theorem the well known
result of Kesten \cite{KE}.
\end{remark}

The interest in the diffusions satisfying the assumptions of Fact
1 stems from the fact that they exhibit subdiffusive behavior;
they are very slow. Compare the $(\log t)^a$ in Fact 1 with the
$t^{1/2}$ for the analogous result for Brownian motion. First
Sinai (\cite{SI}) studied the discrete time analog of diffusion in
Brownian environment, which is random walk in random environment
on $\D{Z}$, and established the analogous to Fact 1 result with
normalizing factor $(\log t)^2$. Shortly after, S. Schumacher
studied the continuous time case.

When $w$ is a stable process having index $a$ and satisfying the
assumptions of Fact 1, $\omega_t/(\log t)^a$ converges in
distribution to $b_1(w)$. The density of $b_1(w)$ has been
computed under various assumptions for $w$. Kesten (\cite{KE})
considered the case where $w$ is a two sided Brownian motion.
Golosov (\cite{GO}), the case where $w(s)=+\infty$ for $s<0$ and
$(w(s))_{s\ge0}$ a Brownian motion, i.e., there is a reflecting
barrier at zero. Tanaka (\cite{TA1}), the case where $w$ is a
symmetric stable process. In the cases where $w$ equals $|B|$ or
$-|B|$, with $B$ two sided Brownian motion having $B_0=0$,
Schumacher's theorem does not apply (in the first case, (iii)
fails; in the second, both (iii) and (iv) fail). However, Tanaka
showed that, in both cases, $\omega(t)/(\log t)^2$ converges in
distribution to a symmetric random variable. For the restriction
of each of the two variables on $[0,+\infty)$ he gave the Laplace
transform. Ours is the first asymmetric case considered.

\section{Preliminaries} \label{prelim}

\subsection{Definition of the process $b$}
For a function $w:\mathbb{R}\to\D{R},\ x>0$, and $y_{0}\in\D{R}$,
we say that \textbf{$w$ admits an $x$-minimum at $y_{0}$} if there
are $\ga,\gb\in\mathbb{R}$ with $\ga<y_{0}<\gb$,
$w(y_{0})=\inf\{w(y):y\in[\ga,\gb]\}$, $w(\ga)\ge w(y_{0})+x$, and
$w(\gb)\ge w(y_{0})+x$. We say that \textbf{$w$ admits an
$x$-maximum at $y_{0}$} if $-w$ admits an $x$-minimum at $y_{0}$.

For convenience, we will call a point where $w$ admits an
$x$-maximum or $x$-minimum, an $x$-maximum or an $x$-minimum
respectively.

We denote by $R_x(w)$ the set of $x$-extrema of $w$ and define
\[
\C{W}_1:=\left\{w\in\C{W}: \begin{array}{c}\text{ For every $x>0$
the set $R_x(w)$ has no accumulation point in $\D{R}$,}  \\
\text{it is unbounded above and below,}\\
\text{and the points of $x$-maxima and $x$-minima
alternate.}\end{array}\right\}.
\]
Thus, for $w\in\C{W}_1$ and $x>0$, we can write
$R_x(w)=\{x_k(w,x):k\in \D{Z}\}$ with $(x_k(w,x))_{k\in\D{Z}}$
strictly increasing, $x_0(w,x)\le0<x_1(w,x)$,
$\lim_{k\to-\infty}x_k(w,x)=-\infty,$ and
$\lim_{k\to\infty}x_k(w,x)=\infty$. Whenever $\PPP(\C{W}_1)=1$,
which will be the case for us, the definition Schumacher gave for
$b$ agrees with the following.

\begin{definition}
The process $b:[0,+\infty)\times\C{W}\to \D{R}$ is defined for
$x>0$ and $w\in \C{W}_1$ as
$$
b_{x}(w):=\begin{cases} x_0(w,x) &  \text{ if $x_0(w,x)$ is an $x$-minimum,} \\
                      x_1(w,x) &  \text{ else,}
        \end{cases}
$$
and $b_x(w)=0$ if $x=0$ or $w\in\C{W}\setminus\C{W}_1$.
\end{definition}

\subsection{Some useful facts}

In the case under consideration, the process $w$ appearing in
\eqref{generator} is a two sided stable process with index
$a\in(1,2]$, no positive jumps, and $w(0)=0$. By two sided stable
we mean that we take two i.i.d stable processes $Y, \widetilde Y$
with paths in $D([0, +\infty)):=\{f\in \D{R}^{[0, +\infty)}: f
\text{ right continuous with left limits}\, \}$, and define $w$ by
$w_s=Y_s$ for $s\ge0$ and $w_s=-\widetilde Y_{(-s)-}$ for $s<0$.
Then $w$ has paths right continuous with left limits. Also it
takes both positive and negative values. In fact for all $t\ge0$,
we have $\PPP(w_t\ge0)=a^{-1}$ (see \cite{BE1}, VIII.1).

Before stating two more properties of $w$, we remind the reader
that for the given L\'evy process, a point $x$ is called regular
for a set $A\subset\D{R}$ if $\PPP_x(\inf\{t>0:w_t\in A\}=0)=1$,
where $\PPP_x$ is the law of the process starting at time 0 from
$x$.

\nid \textbf{Fact 2.}

\begin{enumerate}
\item [(\textsc{i})] 0 is regular for $(0,+\infty)$ and
$(-\infty,0)$.

\item [(\textsc{ii})] $\varliminf_{t\to\pm\infty} w_t=-\infty,
\varlimsup_{t\to\pm\infty} w_t=+\infty$.
\end{enumerate}

The first statement follows from Rogozin's criterion (Proposition
VI.11 in \cite{BE1}) and the fact that
$\int_0^1t^{-1}\PPP(w_t\ge0)dt=\int_0^1t^{-1}\PPP(w_t\le0)dt=+\infty$.
The second statement follows from Theorem VI.12 in \cite{BE1} and
the fact that
$\int_1^{+\infty}t^{-1}\PPP(w_t\ge0)\,dt$=$\int_1^{+\infty}t^{-1}\PPP(w_t\le0)\,dt=+\infty$

The assumptions of Fact \ref{SchumProp} are satisfied. Relation
\eqref{noexplosion} holds because of Fact 2 (\textsc{ii}), and
every process in the set $\{(w(st)/t^{1/a})_{s\in\D{R}}:t>0\}$ has
the same law as $w$. So $U=w$. Then (ii) and (iii) follow from
Lemma \ref{quasicont} (see Section \ref{lemmata}), whose
assumption holds (Fact 2 (\textsc{i}) above), while (i) and (iv)
are clearly true.

Also Fact 2 and Lemma \ref{discrete} (see Section \ref{lemmata})
show that $\PPP(\C{W}_1)=1$, and so the process $b$ is determined
by the definition of the previous subsection.

The absence of positive jumps implies that $\EEE(\exp\{\gl
w_t\})<\infty$ for all $t, \gl>0$ (\cite{BE1}, VII.1). Thus, the
characteristic function of $w_1$ extends to an analytic function
in $\{z\in \D{C}: Im(z)\le 0\}$, and by its form (VIII.1 in
\cite{BE1}) we can see that there is a positive constant $c$ such
that
\[
\EEE(\exp\{\gl w_t\})=\exp\{ct\gl^a\} \text{\qquad  for all $t,\gl
>0$.}
\]
In this work, we assume that $c=1$ as every other case reduces to
this one after a normalization.

\begin{remark}We stick to the spectrally negative case because for $w$
spectrally positive, the process $\tilde w$ defined by $\tilde
w_t=\lim_{s\nearrow -t}w_s$ for $t\in\D{R}$ is spectrally negative
stable with the same index and $b_\cdot(w)=-b_\cdot(\tilde w)$.
\end{remark}

\section{Preparation and proof of the Theorem } If $w$ is under $\PPP$
 a two sided stable process with
$\EEE(\exp\{\gl w_t\})=\exp\{t\gl^a\}$ for all $t,\gl\ge0$ and
some $a\in(1, 2]$, then $\PPP(\C{W}_1)=1$ as was explained in the
previous section. So for $x>0$ let $R_x(w)=\{x_k:k\in\mathbb{Z}\}$
be the set of $x$-extrema for $w$, with $(x_k)_{k\in\mathbb{Z}}$
strictly increasing and $x_0\le0<x_1$.

\begin{lemma}\label{stablerenewal}
The trajectories between consecutive $x$-extrema,
$(w_{x_k+t}-w_{x_k}: t\in[0, x_{k+1}-x_k])$ with $k\in\mathbb{Z}$,
are independent, and the ones corresponding to even non zero $k$
(resp. odd $k$) are identically distributed.
\end{lemma}

 The proof of the lemma is given in Section \ref{lemmata}.
 We call the translation $(w-w(x_k))|[x_k,x_{k+1}]$ of the
trajectory of $w$ between two consecutive $x$-extrema an
\textbf{$x$-slope} (or a slope, when the value of $x$ is clear or
irrelevant), a slope that takes only non-negative values an
\textbf{upward slope}, and a slope taking only non-positive values
a \textbf{downward slope}. We call $(w-w(x_0))|[x_0,x_1]$ the
\textbf{central $x$-slope}.

\begin{remark}As will become clear when we will examine the structure of the
$k$-slopes in Section \ref{hittingtimescomputations}, the upward
and downward $k$-slopes of a two sided Brownian motion (i.e.,
$a=2$), excluding the central one, are identically distributed up
to sign change. When $a\ne2$, since the process has only negative
jumps, the upward $k$-slopes are essentially different from the
downward.
\end{remark}

For any $x$-slope $T:[\alpha,\beta]\to\D{R}$ we call
$l(T):=\beta-\alpha$  the \textbf{length} of $T$. Also we denote
by $\theta(T)$ the slope with domain $[0,\beta-\alpha]$ and values
$\theta(T)(\cdot)=T(\alpha+\cdot)$.

First we determine the distribution functions $F_u, F_d$ of the
lengths $\ul{l}_1, \ol{l}_1$ of an upward and a downward $1$-slope
respectively from the common distributions mentioned in the
preceding lemma. By scaling, this gives the laws for the
$x$-slopes when $x\ne1$. The proof of the following lemma is given
in Section \ref{hittingtimescomputations}.

\begin{lemma} \label{slopelengths}
For all $u>0$,
\begin{align*}
\EEE(e^{-u\underline{l}_1})&=
\left(\Gamma(a+1)\big(E_a^{\prime}(u)+\frac{a}{a-1}\,u\,E_a^{\prime\prime}(u)\big)\right)^{-1},\\
\EEE(e^{-u\overline{l}_1})&=\Gamma(a+1)\Big(E_a^{\prime}(u)+\frac{a}{a-1}\,u\,E_a^{\prime\prime}(u)-\frac{a}{a-1}\,u\,
\frac{\big(E_a^{\prime}(u)\big)^2}{E_a(u)}\Big).
\end{align*}
In particular, the mean values of $\ul{l}_1, \ol{l}_1$ are
\begin{align*}
\EEE(\underline{l}_1)&=\frac{1}{a-1}\frac{\Gamma(a)}{\Gamma(2a-1)}, \\
\EEE(\overline{l}_1)&=\frac{1}{a-1}\Big(\frac{1}{\Gamma(a)}-\frac{\Gamma(a)}{\Gamma(2a-1)}\Big).
\end{align*}
\end{lemma}

\nid And now we are ready to prove our theorem.

\medskip

 \nid \textbf{Proof
of the Theorem:}

For $t\in\D{R}$, let $T_t$ be the 1-slope around $t$. More
precisely, $T_t$ is the slope with domain$(T_t)=[c_t, d_t]$ and
$t\in[c_t, d_t)$. And define $$q:=\inf\{t>0: T_t \text{ downward
slope with domain}(T_t)\subset (0,+\infty)\}.$$ Then
\begin{multline} \label{t1}
\PPP(b_1(w)>x)=\PPP(T_0 \text{ downward slope and } d_0>x)\\
=\PPP(T_t \text{ downward slope and } d_t-t>x)
\end{multline}
 for all $t>0$ because $T_0(-t+\cdot)$ is the same as the
slope around $t$ for $(w_{s-t}-w_{-t}: s\in \D{R})$, and the
latter process has the same law as $w$. Call $C_t$ the event in
the last probability. Then
\begin{equation} \label{t2}
\PPP(C_t)=\PPP(C_t\text{ and } q< t)+\PPP(C_t \text{ and } q\ge
t).
\end{equation}
The second term goes to $0$ as $t\to+\infty$. To work with the
first term,  we consider a sequence of independent random
variables $(\xi_n)_{n\ge1}$, with $\xi_n$ having distribution
function $F_d$ if $n$ is odd, and $F_u$ if $n$ is even. Also let
$(\gz_n)_{n\ge1}$ be a sequence of independent random variables
with $\gz_n\overset{\text{law}}{=}\xi_{n+1}$ for all $n\ge1$.

\nid Finally, define
\begin{align*}
 S_n&:=\xi_1+\cdots+\xi_n, & \tilde
S_n &:=\gz_1+\cdots+\gz_n,\\
 N(t)&:=\inf\{n\ge1:S_n\ge t\}, & \tilde
N(t)&:=\inf\{n\ge1:\tilde S_n\ge t\},\\
B_t&:=S_{N(t)}-t, & \tilde B_t&:=\tilde S_{\tilde N(t)}-t.
\end{align*}
Then Lemma \ref{stablerenewal} implies that
\begin{equation}\label{t3} \PPP(C_t\text{ and } q<
t)=\int_0^t\PPP(N(t-y) \text{ is odd and } B_{t-y}>x)\ dF_q(y),
\end{equation}
where $F_q$ is the distribution function of $q$. We will show that
the limit
$$\lim_{t\to+\infty}\PPP(N(t) \text{ is odd and } B_t>x)$$ exists,
and we will compute its value.

\nid Define $g_1, g_2:[0,+\infty)\to\D{R}$ by
\begin{align*}g_1(t):&=\PPP(N(t)\text{ is odd and } B_t>x),\\
g_2(t):&=\PPP(\tilde N(t)\text{ is even and } \tilde B_t>x)
\end{align*}
for all $t\ge0$. Then $g_1(t)=\PPP(\xi_1>t+x)+\int_0^t
g_2(t-s)\,dF_d(s)$ and $g_2(t)=\int_0^t g_1(t-s)\,dF_u(s)$.
Consequently, $g_1(t)=\PPP(\xi_1>t+x)+g_1\ast(F_u\ast F_d)$, where
$(F_u\ast F_d)(x):=\int_\D{R} F_u(x-y)\,dF_d(y)$ for all
$x\in\D{R}$ is the distribution function of $\xi_1+\xi_2$, and by
the renewal theorem (see \cite{DU}, Chapter 3, statement (4.9))
$$\lim_{t\to\infty}g_1(t)=\frac{1}{\EEE(\ul{l}_1)+\EEE(\ol{l}_1)}\int_0^{+\infty}\PPP(\xi_1>s+x)\,ds=
\frac{1}{\EEE(\ul{l}_1)+\EEE(\ol{l}_1)}\int_x^{+\infty}\PPP(\xi_1>s)\,ds.$$
By \eqref{t1}, \eqref{t2}, \eqref{t3}, this equals
$\PPP(b_1(w)>x)$. Differentiating with respect to $x$ and noting
that $\EEE(\ul{l}_1)+\EEE(\ol{l}_1)=((a-1)\,\Gamma(a))^{-1}$ (by
Lemma \ref{slopelengths}), we find the density of $b_1(w)$ in
$[0,+\infty)$ as stated in the Theorem. The density in
$(-\infty,0]$ is found similarly. \hfill \qedsymbol

\section{Hitting times computations}
\label{hittingtimescomputations} According to Lemma
\ref{stablerenewal}, excluding the central $k$-slope, the images
of all upward (resp. downward) $k$-slopes under the map $\theta$
have the same distribution say $m_{k,u}$ (resp. $m_{k,d}$). In
this section we describe the structure of a $k$-slope picked from
either distribution.

Consider a L\'evy process $(X_t)_{t\ge0}$ starting from 0 and for
which 0 is regular for $(-\infty, 0)$ and $(0,+\infty)$.

\nid For $t>0$ define
\begin{align*}\ul{X}_t&:=\inf\{X_s:s \in [0,t]\},\\
\ol{X}_t&:=\sup\{X_s:s \in [0,t]\}, \intertext{and for $k>0$,}
\ol{\tau}_k&:=\inf\{t>0:\ol{X}_t-X_t\ge k\},\\
\ol{\sigma}_k&:=\sup\{s<\ol{\tau}_k:\ol{X}_s=X_s\},\\
\overline{\beta}_k&:=\ol{X}_{\ol{\tau}_k},\\
\\
\ul{\tau}_k&:=\inf\{t>0:X_t-\ul{X}_t\ge k\},\\
\ul{\gs}_k&:=\sup\{s<\ul{\tau}_k:\ul{X}_s=X_s\},\\
\ul{\gb}_k&:=-\ul{X}_{\ul{\tau}_k}.
\end{align*}

If $\ol{\tau}_k,\ul{\tau}_k<+\infty$ a.s., then $X$ is continuous
at $\ol{\gs}_k, \ul{\gs}_k$ , and splitting the path of $X$ at
$\ol{\gs}_k$ (or $\ul{\gs}_k$) creates two independent pieces (see
lemmata \ref{quasicont} and \ref{decomp} in Section 5).

In the following, we will use the operation of ``gluing together''
functions defined on compact intervals. For two functions $f:[\ga,
\gb]\to\D{R}$, $g:[\gamma, \gd]\to\D{R}$, by gluing $g$ to the
right of $f$ we mean that we define a new function $j:[\ga,
\gb+\gd-\gamma]\to\D{R}$ with
\begin{equation*}
j(x)=
 \begin{cases}
f(x) &\mbox{ for } x\in [\ga, \gb], \\
f(\gb)+g(x-\gb+\gamma)-g(\gamma) &\mbox{ for } x\in [\gb,
\gb+\gd-\gamma].
 \end{cases}
\end{equation*}
Clearly an upward $k$-slope picked from $m_{k,u}$ is obtained by
gluing two independent trajectories with law
$(X_{\ul{\gs}_k+s}-X_{\ul{\gs}_k}:s\in[0,\ul{\tau}_k-\ul{\gs}_k])
,\ (X_s:s\in[0,\ol{\gs}_k])$ in this order, while a downward
$k$-slope picked from $m_{k,d}$ is obtained by gluing two
independent trajectories with law
$(X_{\ol{\gs}_k+s}-X_{\ol{\gs}_k}:s\in[0,\ol{\tau}_k-\ol{\gs}_k])
,\ (X_s:s\in[0,\ul{\gs}_k])$ in this order.

In the remaining of this section, we compute the Laplace
transforms of the distributions of the lengths and heights of
these four kinds of trajectories in the case that $X$ is a L\'evy
process with no positive jumps, and for which
$\ol{\tau}_k,\ul{\tau}_k<+\infty$ a.s. In particular, we exclude
the case where $X$ is the negative of a subordinator. As already
mentioned in section \ref{prelim}, the absence of positive jumps
implies that $\EEE(\exp\{\gl X_t\})<\infty$ for all $\gl>0$ (see
\cite{BE1} VII.1). Let $\psi:[0,+\infty)\to\D{R}$ be defined by
$$ \EEE(\exp\{\gl X_t\})=\exp\{t\psi(\gl)\}
$$
for all $\gl\ge0$. It holds that $\psi$ is convex with $\psi(0)=0,
\psi(+\infty)=+\infty$ (see Chapter VII in \cite{BE1}). Denote by
$\Phi(q)$ the largest root of $\psi(x)=q$. For every $q\ge0$ there
is a continuous function $W^{(q)}:[0,+\infty)\to[0,+\infty)$ with
Laplace transform
\begin{equation}\label{qscale}
\int_0^{+\infty}e^{-\gl x}W^{(q)}(x)\,dx=\frac{1}{\psi(\gl)-q}
\qquad \text{ for all }\gl>\Phi(q).
\end{equation}
The family of functions $\{W^{(q)}: q\ge0\}$ appears in the
solution of the exit problem for $X$. More specifically, if for
$0<x<y$ we define $T:=\inf\{t>0:X_t\notin (0,y)\}$, then (see
\cite{BE3})
\begin{equation}\label{LaplEx}
\EEE_x(e^{-qT}1_{X_T=y})=W^{(q)}(x)/W^{(q)}(y) \text{\qquad for
all $q\ge0$.}
\end{equation}
For every $q\ge0$, we define the function $Z^{(q)}:[0,+\infty)\to
[1,+\infty)$ by
$$Z^{(q)}(x)=1+q\int_0^x W^{(q)}(z)\,dz.$$
\textbf{Note:} In the following, instead of $W^{(0)}, Z^{(0)}$ we
write just $W, Z$.

\medskip

We also introduce a family of processes that is obtained from $X$
by a change of measure. More specifically, since for $c\ge0$ the
process $(\exp\{cX_t-\psi (c) t\})_{t\ge0}$ is a martingale with
mean 1, we can introduce the probability measure $\PPP^c$ for
which
\[
\frac{d\PPP^c}{d\PPP}\Big|_{\C{F}_t}=\exp\{cX_t-\psi (c) t\}
\text{\qquad for all }t\ge0 .
\]
 It is easy to see that $X$ is under $\PPP^c$ a L\'evy process with no positive jumps for
which $0$ is regular for $(-\infty,0)$ and $(0,+\infty)$. Its
Laplace exponent is given by $\psi_c(\gl)= \psi(\gl+c)-\psi(c)$
for $\gl\ge0$. We denote by $\EEE^c$ the expectation with respect
to $\PPP^c$ and
 by $\Phi_c ,W_c^{(q)},Z_c^{(q)}$ the corresponding functions.
As proved in \cite{AKP}, for fixed $c, x\ge0$, the map
$W_c^{(\cdot)}(x)$ can be extended uniquely to an analytic
function in $\D{C}$. The same holds for $Z_c^{(\cdot)}(x)$
obviously. A relation that we will use in the following is
\begin{equation}\label{scalerelation}
W^{(u-\psi(c))}_c(x)=e^{-cx}W^{(u)}(x)
\end{equation}
for all $x,c\ge0, u\in\D{C}$. It is Remark 4 in \cite{AKP}.

\medskip

\nid The main result of this section is the following.
\begin{proposition} \label{stableProp}
Let $X$ be a spectrally negative L\'evy process for which zero is
regular for $(-\infty, 0)$ and $(0,+\infty)$, and $k>0$ such that
$\ol{\tau}_k, \ul{\tau}_k<\infty$ a.s. Then for $u,v\ge0$ it holds
\begin{align}
\EEE(\exp\{-u(\ol{\tau}_k-\ol{\gs}_k)-v(\ol{X}_{\ol{\tau}_k}-X_{\ol{\tau}_k}-k)\})&=e^{vk}
\frac{W(k)}{W'(k)}\Big(Z_v^{(p)}(k)\frac{W_v^{(p)\prime}(k)}{W_v^{(p)}(k)}-pW_v^{(p)}(k)\Big)\label{lfirst},\\
\EEE(\exp\{-u\ol{\gs}_k\}1_{\overline{\beta}_k\le
x})&=\frac{W^{\prime}(k)}{W(k)}\frac{W^{(u)}(k)}{W^{(u)\prime}(k)}\Big(1-e^{-x\frac{
W^{(u)\prime}(k)}{W^{(u)}(k)}}\Big) \label{lsecond}, \\
\EEE(\exp\{-u(\ul{\tau}_k-\ul{\gs}_k)\})&=
\frac{W(k)}{W^{(u)}(k)}\label{lthird},\\
\EEE(\exp\{-u\ul{\gs}_k\})&=\frac{1}{Z^{(u)}(k)}\frac{W^{(u)}(k)}{W(k)}\label{lfourth},
\\
\EEE(\exp\{-u\ul{\beta}_k\})&=\frac{e^{-uk}}{Z_u^{(-\psi(u))}(k)}\label{lfifth},
\end{align}
where $p=u-\psi(v)$.
\end{proposition}
\begin{proof}
 Most of the formulas are contained in the computations in \cite{AKP} and \cite{PI}. We
provide the parts not treated there. \nid Relation~\eqref{lfirst}
is relation (17) in page 223 of \cite{AKP}.

\nid Relation~\eqref{lsecond} is proved by modifying the argument
in the computation of $I_1$ in pages 222, 223 of the same paper.
That is, we integrate up to local time $x$.

\nid For the proof of relation~\eqref{lthird} we will use facts
and the standard notation from excursion theory (see e.g.
\cite{BE1} Chapter IV). Denote by $D[0,+\infty)$ the space of real
valued functions with domain $[0,+\infty)$ which are right
continuous and have left limits everywhere. The set
$$ \C{E}:=\{\gep \in D[0,+\infty) : \exists a\in (0,+\infty] \text{ such that
$\gep(x)\ne0$ for } x\in (0,a)\}$$ is called the space of
excursions. Together with $\C{E}$, we introduce a new point $\gd$
whose use will appear shortly.  Let $(L(t))_{t\ge 0}$ be a local
time process at zero for $Y:=X-\ul{X}$, and define
$L^{-1}(t):=\inf\{s\ge0:L(s)>t\}$ for all $t>0$. The excursion
process $(e_t)_{t>0}$ of $Y$ is given by
$$e_t(s):=
 \begin{cases}
 Y_{L^{-1}(t-)+s} 1_{0\le s\le L^{-1}(t)-L^{-1}(t-)} &\mbox{ if } L^{-1}(t)-L^{-1}(t-)>0, \\
 \gd  &\mbox{ otherwise,}
 \end{cases}
$$
for all $t>0,s\ge0$. It is a Poisson point process with values in
$\C{E}\cup \{\gd\}$; we denote by $n$ its characteristic measure.
For $\gep\in\C{E}$, we call $\ol{\gep}:=\sup_{s\in[0,\gz]}\gep(s)$
the height of $\gep$. Now returning to what we want to compute,
observe that the random variable $\ul{\tau}_k-\ul{\gs}_k$ is the
time that it takes for the first excursion of $Y$ with height at
least $k$ to reach $k$. The law of this excursion is
$n(\gep|\overline{\gep}\ge k)$ (see \cite{BE1}, Chapter 0,
Proposition 2). Thus, the expectation we want is
$$\int_{\mathcal{E}}e^{-u\rho_k}\,dn(\gep|\overline{\gep}\ge
k),$$ where $\rho_k:=\inf\{t\ge 0 : \gep(t)\ge k\}$. For $\theta
\in (0,k]$, we define $\rho_\theta:=\inf\{t\ge 0 :
\gep(t)\ge\theta\}$,
$\mathcal{G}_\theta:=\sigma(\gep(t):t\le\rho_\theta)$, and
$$M_\theta:=e^{-u\rho_\theta}\frac{W^{(u)}(\theta)}{W^{(u)}(k)}\frac{W(k)}{W(\theta)}.$$
The denominator is not zero due to the assumption
$\ul{\tau}_k<\infty$ a.s. and \eqref{LaplEx}. We claim that
$(M_\theta)_{\theta\in(0,k]}$ is a martingale with respect to the
measure $n_k:=n(\cdot|\overline{\gep}\ge k)$ and the filtration
$(\mathcal{G}_\theta)_{\theta\in(0,k]}$. Denote by $\EEE_\mu$ the
expectation with respect to any given measure $\mu$. Observe that
$M_k=e^{-u\rho_k}$ and
\begin{equation}\label{marting}
\EEE_{n_k}(M_k|\mathcal{G}_\theta)=
\EEE_n(e^{-u\rho_k}1_{\rho_k<\infty}|\mathcal{G}_\theta)/\EEE_n(1_{\rho_k<\infty}|\mathcal{G}_\theta).
\end{equation}
Using the Markov property for excursions (see Theorem VI.48.1 in
\cite{RW}), the absence of positive jumps, and \eqref{LaplEx}, we
see that the numerator equals
$$e^{-u\rho_\theta}\EEE_\theta(e^{-uT_k^+}1_{T_k^+<T_0^-})=e^{-u\rho_\theta}W^{(u)}(\theta)/W^{(u)}(k),$$
where $\EEE_\theta$ is the expectation with respect to the law of
$X$ starting from $\theta$, and
\begin{align*}
T_k^+:&=\inf\{t\ge0:X_t\ge k\},\\
 T_0^-:&=\inf\{t\ge0:X_t\le 0\}.
\end{align*}
Now set $u=0$ in the last expression to find the value of the
denominator in \eqref{marting} as
$$\EEE_n(1_{\rho_k<\infty}|\mathcal{G}_\theta)=W(\theta)/W(k).$$
Thus, $\EEE_{n_k}(M_k|\mathcal{G}_\theta)=M_\theta$ proving the
claim. So
$$\EEE_{n_k}(e^{-u\rho_k})=\EEE_{n_k}(M_k)=\EEE_{n_k}(M_\theta)=
\frac{W^{(u)}(\theta)}{W^{(u)}(k)}\frac{W(k)}{W(\theta)}\EEE_{n_k}(e^{-u\rho_\theta})$$
for all $\theta\in(0,k]$. Relation \ref{lthird} will be proved if
we show that the limit of the last quantity as $\theta\downarrow
0$ is $W(k)/W^{(u)}(k)$. Certainly
$\lim_{\theta\downarrow0}\EEE_{n_k}(e^{-u\rho_\theta})=1$ using
the bounded convergence theorem. For the term
$W^{(u)}(\theta)/W(\theta)$ we use relation (9) from \cite{BE3},
which is $W^{(u)}(\theta)=W(\theta)+\sum_{j=1}^{+\infty} u^j
W^{*(j+1)}(\theta)$, and the bound $W^{*(j+1)}(\theta)\le \theta^j
W(\theta)^{j+1}/j!$ (relation (10) in \cite{BE3}). These give that
$\lim_{\theta\downarrow0} W^{(u)}(\theta)/W(\theta)=1$.

 Relation~\eqref{lfourth} follows from \eqref{lthird}, the independence of $\ul{\sigma}_k, \ul{\tau}_k-\ul{\sigma}_k$
 (Lemma \ref{decomp}), and
 the expression for the Laplace transform of $\ul{\tau}_k$ given in Proposition 2 of \cite{PI} as
\begin{equation}\label{Pist}
 \EEE(\exp\{-u\ul{\tau}_k\})=\frac{1}{Z^{(u)}(k)} \text{\qquad for } u\ge0.
\end{equation}
Finally, for relation \eqref{lfifth} we compute
$\EEE(\exp\{-u\ul{\beta}_k\})=\EEE(\exp\{u(k-\ul{\beta}_k)-uk\})=e^{-uk}\EEE(\exp\{uX_{\ul{\tau}_k}\})$.
For $n\in\D{N}$ we have
\begin{multline*}
\EEE(\exp\{uX_{(\ul{\tau}_k\wedge
n)}\})=\EEE(\exp\{uX_{(\ul{\tau}_k\wedge
n)}-\psi(u)(\ul{\tau}_k\wedge n)+\psi(u)(\ul{\tau}_k\wedge
n)\})\\= \EEE^u(\exp\{\psi(u)(\ul{\tau}_k\wedge n)\}).
\end{multline*} Taking $n\to+\infty$ and applying the dominated
convergence theorem in the first quantity (since
$X_{\ul{\tau}_k}\le k$ by the absence of positive jumps) and the
monotone convergence theorem in the last quantity of the last
relation, we obtain
\begin{equation}\label{PistExt}
\EEE(\exp\{uX_{\ul{\tau}_k}\})=
\EEE^u(\exp\{\psi(u)\ul{\tau}_k\}).
\end{equation}
To compute the last expectation, observe that relation
\eqref{Pist} written for the spectrally negative L\'evy process
$(X, \PPP^u)$ is
\begin{equation}\label{PistPreExt}
\EEE^u(\exp\{-v\ul{\tau}_k\})=\frac{1}{Z_u^{(v)}(k)} \text{\qquad
for }v\ge0.
\end{equation}
We will show that this holds for $v=-\psi(u)$ also. So assume that
$\psi(u)>0$. The left hand side is finite for $v=-\psi(u)$ (due to
\eqref{PistExt} and $X_{\ul{\tau}_k}\le k$), which implies that it
can be written as a power series of $-v$ in $D(0, \psi(u))$ with
positive coefficients, continuous in $\overline{D(0, \psi(u))}$.
The denominator of the right-hand side can be extended to an
entire function of $v$ as was mentioned just before this
proposition. By a well known property of analytic functions,
equation \eqref{PistPreExt} will hold for all $v\in\overline{D(0,
\psi(u))}$. Applying it for $v=-\psi(u)$ and combining it with the
equalities before it, we obtain \eqref{lfifth}.
\end{proof}

\nid In the case of our interest, $\psi(u)=u^a$. Consequently,
\begin{align*}
\Phi(u)&=u^{1/a}, &
W(z)&=z^{a-1}/\Gamma(a) ,\\
W^{(u)}(z)&=az^{a-1}E_a^{\prime}(uz^a),  & Z^{(u)}(k)&=E_a(uk^a),
\end{align*}
for $z,u,v\in[0,+\infty)$. The expressions for $W,\ W^{(u)}$ are
found in \cite{BE2}.

\medskip

\nid \textbf{Proof of Lemma \ref{slopelengths}:} Remember the
description of $k$-slopes given in the beginning of this section
just after the definition of the ``gluing together" operation,
where the role of $X$ is played now by $w$. It follows that the
length of an upward 1-slope picked from $m_{1,u}$ equals in
distribution to
$\underline{\tau}_1-\underline{\sigma}_1+\overline{\sigma}_1^*$,
where $\overline{\sigma}_1^*$ is independent of
$\underline{\tau}_1-\underline{\sigma}_1$, and
$\overline{\sigma}_1^*\overset{\text{law}}{=}\overline{\sigma}_1.$
Similarly, the length of a downward 1-slope picked from $m_{1,d}$
equals in distribution to
 $\overline{\tau}_1-\overline{\sigma}_1+\underline{\sigma}_1^*$, with $\ul{\sigma}_1^*$ independent of
$\ol{\tau}_1-\ol{\sigma}_1$, and
$\ul{\sigma}_1^*\overset{\text{law}}{=}\ul{\sigma}_1.$ Using
relations \eqref{lfirst}, \eqref{lsecond}, \eqref{lthird}, and
\eqref{lfourth}, we get the formulas for the Laplace transforms.
For the mean values of $\ul{l}_1, \ol{l}_1$, we compute the
derivatives of their Laplace transforms at zero. To justify the
move of the differentiation under the expectation, we use the
monotone convergence theorem; which applies because the length of
a slope is a positive random variable and the function $(\gl
\mapsto(1-e^{-a\gl})/\gl)$ is nonnegative and decreasing in
$\D{R}$ for any $a>0$. \hfill \qedsymbol


\bigskip

\nid \textbf{Justification of Remark \ref{biasCh1}:} Using our
theorem, the fact that $F_u$ is the distribution function for
$\ul{l}_1$, and
$\EEE(\ul{l}_1)+\EEE(\ul{l}_1)=((a-1)\gG(a))^{-1}$, we obtain
$\PPP(b_1<0)=\EEE(\ul{l}_1)/(\EEE(\ul{l}_1)+\EEE(\ul{l}_1))=e^{-g(a)}$,
where
$$g(a):=\log\left(\frac{\EEE(\ol{l}_1)}{\EEE(\ul{l}_1)}+1\right).$$\
By Lemma \ref{slopelengths},
\begin{equation}\label{geq}
g(a)=\log\frac{\gG(2a-1)}{\gG^2(a)}, \end{equation} and using the
formula (see \cite{ER}, \S 1.9, relation (1))
$$\log\gG(z)=\int_0^{+\infty}\left[(z-1)-\frac{1-e^{-(z-1)t}}{1-e^{-t}}\right]
\frac{e^{-t}}{t}\,dt \text{\quad for Re$(z)>0,$}$$ we get
\begin{equation}\label{ERD}
g(a)=\int_0^{+\infty}\frac{e^{-t}}{t(1-e^{-t})}(1-e^{-(a-1)t})^2\,dt.
\end{equation}
$g$ was defined above only for $a\in(1,2]$, but \eqref{geq}
extends it to a differentiable function in $(1/2,+\infty)$. It
follows from \eqref{ERD} that $g$ is a strictly increasing
function of $a$ in $[1, +\infty)$, and \eqref{geq} shows that
$g(1)=0,\ g(2)=\log2$. \hfill \qedsymbol

\section{Some lemmata} \label{lemmata}

\nid In this section we prove some auxiliary results that we used
above.

\smallskip

\nid \textbf{Proof of Lemma \ref{stablerenewal}:} Let
$\ol{\tau}_{x,+}$, $\ol{\gs}_{x,+}$, $\ul{\tau}_{x,+}$,
$\ul{\gs}_{x,+}$ be defined as in the beginning of Section
\ref{hittingtimescomputations} with $w$ having the role of $X$ and
$k=x$. Similarly for $\ol{\tau}_{x,-}$, $\ol{\gs}_{x,-}$,
$\ul{\tau}_{x,-}$, $\ul{\gs}_{x,-}$ with $ \hat w(s)=w(-s)\text{
for } s\ge0$ having the role of $X$ and $k=x$. There are four
possible cases for the ordering of the pairs $\{\ol{\tau}_{x,+},
\ul{\tau}_{x,+}\}$, $\{\ol{\tau}_{x,-}, \ul{\tau}_{x,-}\}$. We
treat only two of them, the other being similar.

First assume that $\ol{\tau}_{x,+}<\ul{\tau}_{x,+}$ and
$\ol{\tau}_{x,-}<\ul{\tau}_{x,-}$. Then $\ul{\gs}_{x,+}$ is a
point of $x$-minimum for $w$, and the path of
$w-w(\ul{\gs}_{x,+})$ after time $\ul{\gs}_{x,+}$ is independent
of the past by Lemma \ref{decomp}.  Similarly, in the negative
semi-axis, $-\ul{\gs}_{x,-}$ is a point of $x$-minimum for $w$ and
breaks the path of $w$ into two independent pieces. Between
$-\ul{\gs}_{x,-}$, $\ul{\gs}_{x,+}$, there is exactly one more
$x$-extremum. It is an $x$-maximum and it is the point in
$\{-\ol{\gs}_{x,-}, \ol{\gs}_{x,+}\}$ where $w$ has greater value.
Say it is $-\ol{\gs}_{x,-}$. Then in the notation of the lemma, we
have $x_{-1}=-\ul{\gs}_{x,-}, x_0=-\ol{\gs}_{x,-}$, and
$x_1=\ul{\gs}_{x,+}$. The points $x_{-1}, x_0, x_1$ cut the path
of $w$ into four independent pieces. This, combined with the fact
that $(w_s)_{s\le0}\overset{\text{law}}{=}(-w_{(-s)-})_{s\le0}$
and time reversal (Lemma II.2 in \cite{BE1}), shows that all
upward $x$-slopes, including $w-w(x_{-1})|[x_{-1}, x_0]$, are
obtained by gluing two trajectories with law
$(w_{s+\ul{\gs}_{x,+}}-w_{\ul{\gs}_{x,+}}:s\in[0,\ul{\tau}_{x,+}-\ul{\gs}_{x,+}])
,\ (w_s:s\in[0,\ol{\gs}_{x,+}])$ in this order in this order,
 while all downward $x$-slopes, excluding $w-w(x_0)|[x_0, x_{-1}]$,  are obtained by gluing two trajectories with law
 $(w_{s+\ol{\gs}_{x,+}}-w_{\ol{\gs}_{x,
+}}:s\in[0,\ol{\tau}_{x,+}-\ol{\gs}_{x, +}]) ,\
(w_s:s\in[0,\ul{\gs}_{x, +}])$. This description accounts for all
$x$-slopes in the path decomposition of $w$.

Now assume that $\ol{\tau}_{x,+}<\ul{\tau}_{x,+}$ and
$\ol{\tau}_{x,-}>\ul{\tau}_{x,-}$. Then $-\ol{\gs}_{x,-}$ is a
point of $x$-maximum for $w$, and $\ul{\gs}_{x,+}$ is a point of
$x$-minimum for $w$. If $w(\ol{\gs}_{x,+})-w(-\ul{\gs}_{x,-})< x$,
then $x_0=-\ol{\gs}_{x,-}, x_1=\ul{\gs}_{x,+}$. If
$w(\ol{\gs}_{x,+})-w(-\ul{\gs}_{x,-})\ge x$, then
$x_0=-\ul{\gs}_{x,-}, x_1=\ol{\gs}_{x,+}$. As in the previous
case, we get the desired description for the decomposition of the
path of $w$ into $x$-slopes. \hfill \qedsymbol

\medskip

For a function $f:[0, +\infty)\to\D{R}$ with only jump
discontinuities and $x_0>0$, we say that $f$ has a left local
maximum (resp. minimum) at $x_0$ if there is an $\gep>0$ such that
$f(x)\le f(x_0-)$ (resp. $f(x)\ge f(x_0-)$) for all
$x\in(x_0-\gep, x_0)$. Similarly for a right local maximum and
minimum.

\begin{lemma}\label{quasicont}
Let X be a L\'evy process such that $0$ is regular for
$(0,+\infty)$ and $(-\infty,0)$. With probability one
\begin{enumerate} \item[(i)] $X$ is continuous at every one sided
local extremum. \item[(ii)] In no two local minima (resp. maxima)
X has the same value.
\end{enumerate}
\end{lemma}
\begin{proof}
$(i).$ It is enough to consider the case of a local one sided
maximum (the case of one sided local minimum follows by applying
the present case to the process $-X$). Consider the set of times
where the process jumps upwards or downwards by at least $1/n$; it
is a countable subset of $(0,+\infty)$ with no accumulation point.
Apply the strong Markov property to each of these times. Since $0$
is regular for $(0,+\infty)$, none of these can be a point of a
right local maximum. Using time reversal (Lemma II.2 in
\cite{BE1}) and the fact that $0$ is regular for $(-\infty,0)$, we
exclude the existence of left local maxima.

$(ii).$ This holds for any L\'evy process that is not compound
Poisson (Proposition VI.4 in \cite{BE1}). It is a simple
application of the strong Markov property.
\end{proof}

\nid In the next lemma, we use the notation introduced in the
beginning of Section \ref{hittingtimescomputations}.
\begin{lemma}\label{decomp}
Let $(X_t)_{t\ge0}$ be a L\'evy process starting from zero such
that $0$ is regular for $(-\infty, 0)$ and $(0, +\infty)$, and
$\varliminf_{t\to\infty} X_t=-\infty$,
$\varlimsup_{t\to\infty}X_t=+\infty$. With probability one:
\begin{enumerate} \item [(i)] The two trajectories
$(X_t:t\in[0,\ol{\gs}_k])$ and
$(X_{\ol{\gs}_k}-X_{\ol{\gs}_k+t}:t\in
[0,\ol{\tau}_k-\ol{\gs}_k])$ are independent. \item[(ii)] The two
trajectories $(X_t:t\in[0,\ul{\gs}_k])$ and
$(X_{\ul{\gs}_k}-X_{\ul{\gs}_k+t}:t\in
[0,\ul{\tau}_k-\ul{\gs}_k])$ are independent.
\end{enumerate}
\end{lemma}
\nid This follows from the discussion in Section 4 of \cite{GP}.
So we don't give its proof. For the next statement, recall that
the set $\C{W}_1$ was defined in Section \ref{intro}.

\begin{lemma}\label{discrete}
If $(w_t)_{t\in\D{R}}$ is a RCLL version of a L\'evy process
starting from zero such that $0$ is regular for $(-\infty, 0)$ and
$(0, +\infty)$, $\varliminf_{t\to\pm\infty} w_t=-\infty$, and
$\varlimsup_{t\to\pm\infty} w_t=+\infty$, then $\PPP(\C{W}_1)=1$.
\end{lemma}
\begin{proof}
First we prove that for fixed $x>0$, the set $$C_x:= \{w\in
C(\D{R}): \text{ $R_x(w)$ has the properties appearing in the
definition of $\C{W}_1$}\}$$ has $\PPP(C_x)=1$. To see this,
observe that for $z$ a point of $x$-minimum and
$\ga_z=\sup\{\ga<z:w(\ga)\ge w(z)+x\}$,
$\gb_z=\inf\{\gb>z:w(\gb)\ge w(z)+x\}$ it holds that
$\ga_z<z<\gb_z$ (because $w$ is continuous at $z$ by Lemma
\ref{quasicont}(i)) and there is no other $x$-minimum in $(\ga_z,
\gb_z)$. Indeed, if $\tilde z$ is an $x$-minimum in $(\ga_z,z)$,
then in case $\gb_{\tilde z}>z$ we get that $w$ takes the same
value in two local minima, while in case $\gb_{\tilde z}<z$ we get
$w(\gb_{\tilde z})\ge w(z)+x$. The first case is excluded by Lemma
\ref{quasicont}(ii), and the second contradicts the definition of
$\ga_z$. If $\tilde z$ is an $x$-minimum in $(z,\gb_z)$, then in
case $\ga_{\tilde z}<z$ we get that $w$ takes the same value in
two local minima, while in case $\ga_{\tilde z}>z$ we get
$w(\ga_{\tilde z}-)\ge w(z)+x$. The first case is excluded by
Lemma \ref{quasicont}(ii), and the second contradicts the
definition of $\gb_z$.

Assume that there is a
 strictly monotone, say increasing, sequence $(z_n)_{n\ge1}$ of $x$-minima converging to
$z_\infty\in\D{R}$. Then by the above observation we get
$\varlimsup_{y, \tilde y\nearrow z_\infty}(w(y)-w(\tilde y))\ge x$
implying that $w$ cannot have left limit at $z_\infty$. A
contradiction with the fact that $w$ is RCLL. Similarly if
$(z_n)_{n\ge1}$ is decreasing. So in a set of $w$'s in $\C{W}$
with probability 1, it holds that the set of $x$-minima of $w$ has
no accumulation point. The same holds for the set of $x$-maxima,
and as a result also for $R_x(w)$. Since
$\varliminf_{|t|\to\infty} w_t=-\infty$ and
$\varlimsup_{|t|\to\infty} w_t=+\infty$, it follows that
$\PPP(R_x(w) \text{ is unbounded above and below} )=1$. Now
between two consecutive $x$-maxima (resp. minima) there is exactly
one $x$-minimum (resp. $x$-maximum). Indeed, take $z_1< z_2$ two
consecutive $x$-minima, and call $s_0$ the unique point where $w$
attains its maximum in $[z_1,z_2]$.  Then $w(s_0)\ge \max\{w(x_1),
w(z_2)\}+x$. Because if $w(s_0)< w(x_1)+x$, then $\gb_{z_1}>z_2$,
while if $w(s_0)< w(z_2)+x$, then $\ga_{z_2}<z_1$, and both
$\gb_{z_1}>z_2$, $\ga_{z_2}<z_1$ are false as was shown above. So
$s_0$ is an $x$-maximum. There is no other x-maximum in $[z_1,
z_2]$ because then we would find an $x$-minimum in $(z_1, z_2)$,
which cannot happen since $z_1, z_2$ are consecutive $x$-minima.
Consequently, $\PPP(C_x)=1$.

Finally, note that for all $n\in\D{N}\setminus \{0\}$ we have
$R_n(w)\subset R_x(w)\subset R_{1/n}(w)$ for $x\in [1/n, n]$, from
which it follows that $\C{W}_1=\cap_{x\in(0, +\infty)}
C_x=\bigcap_{n\in\D{N}\setminus\{0\}}(C_n\cap C_{1/n})$. Thus,
$\PPP(\C{W}_1)=1$.
\end{proof}

\textbf{Acknowledgments.} My advisor, Amir Dembo, and the
anonymous referee pointed out to me some errors in previous drafts
of this paper. I am grateful to both of them. I also thank the
referee for comments that greatly improved the structure and
exposition of the paper.

\bibliographystyle{plain}

\end{document}